\documentstyle[12pt]{article}
\begin{document}
\newpage
\pagestyle{empty}
\setcounter{page}{0}
%
\def\CPbar{\hbox{{\rm CP}\hskip-1.80em{/}}}
\def\hfl#1#2{\smash{\mathop{\hbox to 12mm{\rightarrowfill}}
\limits^{\scriptstyle#1}_{\scriptstyle#2}}}
\def\vfl#1#2{\llap{$\scriptstyle #1$}\left\downarrow
\vbox to 6mm{}\right.\rlap{$\scriptstyle#2$}}
\def\ihfl#1#2{\smash{\mathop{\hbox to 12mm{\nearrowfill}}
\limits^{\scriptstyle#1}_{\scriptstyle#2}}}
\def\ibfl#1#2{\smash{\mathop{\hbox to 12mm{\searrowfill}}
\limits^{\scriptstyle#1}_{\scriptstyle#2}}}
\def\diagram#1{\def\normalbaselines{\baselineskip=0pt
\lineskip=10pt\lineskiplimit=1pt}  \matrix{#1}}
\def\adots{\mathinner{\mkern2mu\raise1pt\hbox{.}
\mkern3mu\raise4pt\hbox{.}\mkern1mu\raise7pt\hbox{.}}}

%
%
%
\def\CC{{\Bbb C}}
\def\NN{{\Bbb N}}
\def\QQ{{\Bbb Q}}
\def\RR{{\Bbb R}}
\def\ZZ{{\Bbb Z}}
\def\cA{{\cal A}}          \def\cB{{\cal B}}          \def\cC{{\cal C}}
\def\cD{{\cal D}}          \def\cE{{\cal E}}          \def\cF{{\cal F}}
\def\cG{{\cal G}}          \def\cH{{\cal H}}          \def\cI{{\cal I}}
\def\cJ{{\cal J}}          \def\cK{{\cal K}}          \def\cL{{\cal L}} 
\def\cM{{\cal M}}          \def\cN{{\cal N}}          \def\cO{{\cal O}}
\def\cP{{\cal P}}          \def\cQ{{\cal Q}}          \def\cR{{\cal R}} 
\def\cS{{\cal S}}          \def\cT{{\cal T}}          \def\cU{{\cal U}}
\def\cV{{\cal V}}          \def\cW{{\cal W}}          \def\cX{{\cal X}}
\def\cY{{\cal Y}}          \def\cZ{{\cal Z}}
\def\qed{\hfill \rule{5pt}{5pt}}
\newtheorem{lemma}{Lemma}
\newtheorem{prop}{Proposition}
\newtheorem{theo}{Theorem}
\newenvironment{result}{\vspace{.2cm} \em}{\vspace{.2cm}}

\rightline{CPTH-9612484}  
\rightline{q-alg/9612028} 
\rightline{December 96} 

\vfill
\vfill

\begin{center}

{\LARGE {

\bf {\sf 
Combined ($q$,$h$)-Deformation as a Nonlinear Map on ${\cal U}_q(sl(2))$
}}}\\[2cm]

\smallskip 

{\large B. Abdesselam$^{\dagger,}$\footnote{boucif@orphee.polytechnique.fr}, 
A. Chakrabarti$^{\dagger,}$\footnote{chakra@orphee.polytechnique.fr} and 
R. Chakrabarti$^{\ddagger}$}

\smallskip 

\smallskip 

\smallskip

{\em  $^{\dagger,}$\footnote{Laboratoire Propre du CNRS UPR A.0014}Centre de Physique 
Th\'eorique, Ecole Polytechnique, \\
91128 Palaiseau Cedex, France.}

\smallskip 

\smallskip 

{\em $^{\ddagger}$Department of Theoretical Physics, University of Madras, 
Guindy campus,  \\
Madras-600025, India.}

\end{center}

\vfill

\begin{abstract}
The generators $({\cal J}_{\pm},{\cal J}_0)$ of the algebra 
${\cal U}_q(sl(2))$ is our starting
point. An invertible nonlinear map involving, apart from $q$, a second
arbitrary complex parameter $h$, defines a triplet $({\hat X},{\hat Y},
{\hat H})$. The
latter set forms a closed algebra under commutation relations. The
nonlinear algebra ${\cal U}_{q,h}(sl(2))$, thus generated, has two different
limits. For $q \rightarrow 1$, the Jordanian $h$-deformation
${\cal U}_{h}(sl(2))$ is obtained. For $h \rightarrow 0$, the $q$-deformed
algebra ${\cal U}_{q}(sl(2))$ is reproduced. From the nonlinear map,
the irreducible representations of the doubly-deformed algebra
${\cal U}_{q,h}(sl(2))$ may be directly and explicitly obtained form the known
representations of the algebra ${\cal U}_q(sl(2))$. Here we consider only 
generic values of $q$.

\end{abstract}

\vfill
\vfill

\pagestyle{plain}

\newpage
\section {Introduction}

The standard $q$-deformation of the $sl(2)$ algebra is familiar
enough. Here we treat this as known. The Jordanian $h$-deformation
${\cal U}_h(sl(2))$ is relatively recent. Its various aspects are being
studied [3-9] intensively. In a recent work \cite{ACC96} the present
authors established a nonlinear invertible map between the generators
of ${\cal U}_{h}(sl(2))$ and the classical $sl(2)$ generators. The
principal interest of the map in \cite{ACC96} was that it simply and
immediately provided a construction of the irreducible representations
of the algebra ${\cal U}_h(sl(2))$. Along similar lines, a map was
established in \cite{ACCH96} relating the $h$-deformed 3-dimensional Euclidean
algebra ${\cal U}_h(e(3))$ and its classical partner $e(3)$.

Here, analogously, we consider an appropriate generalization of the
map in \cite{ACC96} . Instead of the generators of the classical algebra
$sl(2)$, we take the generators of the standard $q$-deformed algebra
${\cal U}_q(SL(2))$ as our starting point~; and, via a nonlinear map,
introduce a second deformation parameter $h$. The generators
$({\hat X},{\hat Y},{\hat H})$, thus achieved, form a closed set upon 
commutation. We refer
to this doubly deformed nonlinear algebra as ${\cal U}_{q,h}(sl(2))$. Such an
algebra may be of interest for various reasons, leading to
applications \cite{R91,ABC}. The algebra ${\cal U}_{q,h}(sl(2))$ has two distinct
limits :  ${\cal U}_{q,h}(sl(2)) \rightarrow {\cal U}_h(sl(2))$ as $q \rightarrow
1$; and ${\cal U}_{q,h}(sl(2)) \rightarrow {\cal U}_q(sl(2))$ as $h \rightarrow
0$. When both limits are imposed the classical algebra $sl(2)$ is, of
course, reproduced. Parallel to our construction \cite{ACC96}, here also the
known representations of ${\cal U}_q(sl(2))$ immediately yield, through the
map, the representations of ${\cal U}_{q,h}(sl(2))$. Here we consider only
the generic values of $q$. Special features \cite{RA89,AC91} for $q$ being a
root of unity will be explored elsewhere. The map to be presented
below provides two induced coalgebraic structures for the doubly
deformed algebra ${\cal U}_{q,h}(sl(2))$. These induced structures reflect the
underlying coalgebraic properties of ${\cal U}_q(sl(2))$ and ${\cal U}_h(sl(2))$
respectively.

\section {The map}

We first briefly review the map \cite{ACC96} between the ${\cal U}_h(sl(2))$ and
$sl(2)$ generators. The complications due to the presence of another
deformation parameter $q (\not= 1)$ may then be easier to grasp. The
generators $(J_{\pm},J_0)$ of the algebra $sl(2)$ satisfy
\begin{eqnarray}
[J_0,J_{\pm}] = \pm J_{\pm} ,\qquad\qquad  [J_+,J_-] = 2J_0
\end{eqnarray}
Introducing an arbitrary complex parameter $h$, we now define
\begin{eqnarray}
&&X = \frac{2}{h} \hbox{arctanh}\frac{hJ_+}{2} = \frac{2}{h} \sum_{n \geq 0}
\frac{(\frac{h}{2} J_+)^{2n+1}}{(2n+1)},\nonumber\\
&&Y = (1 - (\frac{h}{2} \ J_+)^2)^{1/2} J_- (1 - (\frac{h}{2}
J_+)^2)^{1/2},\nonumber\\
&& H = J_0 .
\end{eqnarray}

It may be shown \cite{ACC96} that the triplet $(X,Y,H)$ satisfies the
commutation relations \cite{O92} corresponding to the generators of the
${\cal U}_h(sl(2))$ algebra:
\begin{eqnarray}
&&[H,X] = \frac{1}{h}\sinh hX \\
&&[H,Y] = -\frac{1}{2} \biggl(Y(\cosh hX) +  (\cosh
hX)Y \biggr)\\
&&[X,Y] = 2H .
\end{eqnarray}

The Casimir operator for the algebra ${\cal U}_h(sl(2))$, then, via the map,
reduces to its expression corresponding to the classical algebra
$sl(2)$:
\begin{eqnarray} 
&&C=\frac{1}{2h} [(\sinh hX)Y+Y(\sinh hX)] + \frac{1}{4}(\sinh\ hX)^2
+H^2\nonumber\\
&&\phantom{C}= \frac{1}{2}(J_+ J_- + J_- J_+) + J_0^2
\end{eqnarray}
As demonstrated in \cite{ACC96}, the map (2) now explicitly 
furnishes the
irreducible representations of the ${\cal U}_h(sl(2))$ algebra.

We now proceed to obtain a $q$-deformed generalization of the map
(2). The basic elements of our construction are the generators of
the $q$-deformed algebra ${\cal U}_q(sl(2))$. These generators obey the
following commutation rules
\begin{equation}
[{\cal J}_0,{\cal J}_{\pm}] = \pm {\cal J}_{\pm} ,\qquad\qquad
 [{\cal J}_+,{\cal J}_-] = [2{\cal J}_0],
\end{equation}
where
$$
[x] = \frac{q^x - q^{-x}}{q - q^{-1}}.
$$

The coproduct structure for the generators reads 
\begin{equation}
\Delta({\cal J}_{\pm}) = {\cal J}_{\pm} \otimes q^{{\cal J}_0} + 
q^{-{\cal J}_0} \otimes{\cal J}_{\pm}, \qquad \qquad
\Delta({\cal J}_0) = {\cal J}_0 \otimes 1 + 1 \otimes {\cal J}_0  .
\end{equation}

For later use, we here obtain the following identity
\begin{equation}
[{\cal J}_+^p,{\cal J}_-] = \frac{[p]}{q-q^{-1}} (q^{{\cal J}_0} 
{\cal J}_+^{p-1} q^{{\cal J}_0} -
q^{-{\cal J}_0} {\cal J}_+^{p-1} q^{-{\cal J}_0}) .
\end{equation}
In order to get a closed algebra, we choose the following
$q$-generalization of the map (2)
\begin{eqnarray}
&&\frac{h{\hat X}}{2} = \sum_{n \geq 0} \alpha_n (\frac{h}{2} 
{\cal J}_+)^{2n+1},\\
&&{\hat Y} = (1 - (\frac{h}{2} {\cal J}_+)^2)^{1/2} {\cal J}_- 
(1 - (\frac{h}{2}
{\cal J}_+)^2)^{1/2},\\
&&{\hat H} = {\cal J}_0
\end{eqnarray}
where
\begin{equation}
\alpha_n=\frac{1}{[2n+1]} P_n(\xi) ,\qquad\qquad \xi = \frac{q^2 + q^{-2}}{2}
\end{equation}
The purpose of introducing the Legendre polynomials $P_n(\xi)$, will 
become evident as we proceed. It plays a crucial role. We note that in (10) 
(in 
space of the standard finite, $(2j+1)$ dimensional representations) the 
convergence of the
series is not a problem as ${\cal J}_+$ is a nilpotent operator. The series 
in (10) may be
looked as a particular $q$-generalization of the series for the
$\hbox{arctanh}$ function. The maps for the operators ${\hat Y},{\hat H}$ in 
(11,12) read
formally the same as in (2), which may now be looked as $q = 1$
limit of the present maps. In order to invert the map, we
postulate a series
\begin{equation}
{\cal J}_+ = \frac{2}{h} \sum_{n \geq 0} \beta_n (\frac{h}{2}{\hat X})^{2n+1} .
\end{equation}

Substituting (14) in (10) and comparing terms 
of identical powers on both sides, we obtain
\begin{equation}
\beta_0 = 1 ,\qquad\qquad \beta_n = - \sum_{m=1}^n \alpha_m Z_{n,m}
\hspace{0.5cm} for\ n \geq 1
\end{equation}
where
\begin{equation}
Z_{n,n} = 1\ ,\ Z_{n,m} = \sum_{partitions} \zeta_{m,\{\nu_p\}}
\prod_{p=1}^{n-m} \beta_b^{\nu_p} .
\end{equation}
The sum over ``partitions'' in the rhs of (16) maintains the
combinatorial properties.
\begin{equation}
\sum_{p=1}^{n-m} p\;\nu_p = n-m,\qquad\qquad\sum_{p=1}^{n-m} \nu_p \leq (2m+1)
\end{equation}
The symmetry factor $\zeta_{m,\{\nu_p\}}$ in the rhs of
(16)
reads
\begin{equation}
\zeta_{m,\{\nu_p\}} = \frac{(2m + 1)~!}{(2m+1
  -\sum_{p=1}^{n-m}\nu_p)! \prod_{b=1}^{n-m}(\nu_p)!}
\end{equation}

Thus, the inversion of the series (10) may be looked as a
combinatorial problem ; and, an 
arbitrary coefficient $\beta_n$ in the series (14) 
may now be found
recursively. Again, owing to the nilpotency of ${\hat X}$, the convergence of
the rhs in (14) readily follows. The first few coefficients are
listed below :
\begin{eqnarray}
&&\beta_1 = -\alpha_1\ ,\ \beta_2 = -\alpha_2 + 3\alpha_1^2\ ,\ \beta_3
= -\alpha_3 + 8\alpha_2\alpha_1 - 12\alpha_1^3 \nonumber\\
&&\beta_4 = -\alpha_4 + 10\alpha_3\alpha_1 + 5\alpha_2^2 -
55\alpha_2\alpha_1^2 + 55\alpha_1^4 ,\nonumber\\
&&\beta_5 = -\alpha_5 + 12\alpha_4\alpha_1 + 12\alpha_3\alpha_2 -
78\alpha_3\alpha_1^2 - 78\alpha_2^2\alpha_1 + 364\alpha_2\alpha_1^3 -
273\alpha_1^5, \nonumber\\
&&\beta_6 = -\alpha_6 + 14\alpha_5\alpha_1 + 14\alpha_4\alpha_2 -
105\alpha_4\alpha_1^2 + 7\alpha_3^2 - 210\alpha_3\alpha_2\alpha_1\nonumber\\
&&\phantom{\beta_6 =  } + 560\alpha_3\alpha_1^3 - 35\alpha_2^3 + 840 \alpha_2^2\alpha_1^2 -
2380\alpha_2\alpha_1^4 + 1428\alpha_1^6 .
\end{eqnarray}

In the limit $q \rightarrow 1$ we have $\alpha_n \rightarrow
\frac{1}{2n+1}$ ; and, consequently, it follows from the successive
terms in (15) that
\begin{equation}
\beta_n \rightarrow \frac{2^{2n}(2^{2n}-1)}{(2n)!} B_{2n},
\end{equation}
where $B_{2n}$ are the Bernoulli numbers. Thus, the series
(14)
may be viewed as a particular $q$-generalization of the $\hbox{tanh}$
function. In this sense, the series (10) and (14) may be expressed 
respectively as
\begin{eqnarray}
&& \frac{h{\hat X}}{2} = \hbox{arctanh}_q \left(\frac{h{\cal J}_+}{2}\right),\\
&&\frac{h{\cal J}_+}{2} = \hbox{tanh}_q \left(\frac{h {\hat X}}{2}\right)
\end{eqnarray}
It is to be emphasized that we define the $q$-functions $\hbox{arctanh}_q$
and $\hbox{tanh}_q$ only by their series expressions in (10) and
(14)
respectively. (Other definitions can and do exist. But the foregoing ones are necessary for our purpose.)

We are now in a position to demonstrate the commutation relations for
the generators $({\hat X},{\hat Y},{\hat H})$ defined by the map (10, 11, 12). 
To this end, we obtain by employing the identity (9)
\begin{eqnarray}
[{\hat X},{\cal J}_-] = \frac{1}{q-q^{-1}} \left[q^{{\cal J}_0} F \left(\frac{h{\cal J}_+}{2}\right)
q^{{\cal J}_0} - q^{-{\cal J}_0} F \left(\frac{h{\cal J}_+}{2}\right) q^{-{\cal J}_0}\right]
\end{eqnarray}
where
\begin{eqnarray}
&& F(x) = \sum_{n \geq 0} P_n(\xi) x^{2n}\nonumber\\
&&\phantom{F(x)}= (1 - 2\xi x^2 + x^4)^{-1/2}\nonumber\\
&&\phantom{F(x)}= (1 - (qx)^2)^{-1/2} (1 - (q^{-1}x)^2)^{-1/2}
\end{eqnarray}
Use of (23) and the identity
\begin{equation}
q^{\pm {\cal J}_0} f({\cal J}_{\pm}) q^{\mp {\cal J}_0} = 
f(q^{\pm 1} {\cal J}_+)
\end{equation}
leads to
\begin{equation}
[{\hat X},{\hat Y}] = [2{\hat H}] .
\end{equation}

This is the direct $q$-deformation of the commutator (5). We may 
point out that the specific choice of the coefficients $\alpha_n$
in (13) is instrumental in ensuring (26), while a 
simple prescription for ${\hat Y}$ is maintained. The other commutators for the
triplet ${\hat X},{\hat Y},{\hat H}$ may also be determined :
\begin{eqnarray}
&&[{\hat H},{\hat X}]=\frac{2}{h}\sum_{n\geq 0} \frac{(2n+1)}{[2n+1]} P_n(\xi)
\left(\frac{h}{2} {\cal J}_+ \right)^{2n+1} ,\\
&&
[{\hat H},{\hat Y}] = -\frac{1}{2} \left(
\frac{1+(\frac{h}{2} {\cal J}_+)^2}{1-(\frac{h}{2} {\cal J}_+)^2} 
{\hat Y}+{\hat Y}\frac{1+(\frac{h}{2}{\cal J}_+)^2}{1-(\frac{h}{2}
{\cal J}_+)^2}\right) 
\end{eqnarray} 

In the rhs of (27) and (28) the operator
${\cal J}_+$ stands for the series expression (14); and this
leads to the closed algebra for the set $({\hat X},{\hat Y},{\hat H})$. It has 
been verified explicitly that the rhs of (26), (27) and (28) are indeed
compatible with the Jacobi identity.  

We now describe an alternate approach to the inversion of the series
(10). Starting with a variant of (27), we
obtain a closed expression for ${\cal J}_+$ in terms of ${\hat X}$ and 
$q^{\pm {\hat H}}$. To this end, we define
\begin{eqnarray}
 && u = \frac{1}{q-q^{-1}} (q^{\hat H} \frac {h}{2} {\hat X} 
q^{-{\hat H}} - q^{-{\hat H}}\frac {h}{2} {\hat X} q^{\hat H}),\\
&& v = \frac {h}{2} {\cal J}_+
\end{eqnarray}
In the $q \rightarrow 1$ limit, we have
\begin{eqnarray}
&& u \rightarrow \frac {h}{2} [H,X] \\
&&\phantom{u \rightarrow } = \frac {1}{2} \sinh hX
\end{eqnarray}
and $v \rightarrow \frac {h}{2} J_+ $. The equality in (31) follows from 
(3). In this limit the two operators
$u$ and $v$ interrelate, via the map (2), as follows
\begin{eqnarray}
&&v = \hbox{tanh}\ \left(\frac {1}{2} \sinh^{-1} 2u\right)\\
&& \phantom{v}= -\frac{1}{2u} + \sqrt{\frac{1}{4u^2} + 1}.
\end{eqnarray}
When $q \not= 1$, the "deformed" relationship between $u$ and $v$ may
be observed as follows. Equation (7,10,12) yield
\begin{eqnarray}
&&u = \sum_{n \geq 0} P_n(\xi) v^{2n+1}\\
&&\phantom{v} = \frac{v}{\sqrt{(1 - (qv)^2) (1 - (q^{-1}v)^2)}}.
\end{eqnarray}
Solving $v$ in terms of $u$, we get
\begin{equation}
v^2 = \xi + \frac{1}{2u^2} \pm \sqrt{\left(\xi + \frac{1}{2u^2}\right)^2 - 1}
\end{equation}
Noticing that $u \rightarrow 0$ as $v \rightarrow 0$, indicates that
we choose
\begin{equation}
v^2 = \xi + \frac{1}{2u^2} - \sqrt{\left(\xi + \frac{1}{2u^2}\right)^2
  - 1}.
\end{equation}
In the limit $q \rightarrow 1$, the expression (33) is
obtained. Equation (38), therefore, embodies the
$q$-generalized closed form expression of the operator ${\cal J}_+$ in terms
of the operators $({\hat H},{\hat X})$.

To sum up, starting from the standard ${\cal U}_q(sl(2))$ algebra, we
constructed a triplet $({\hat X},{\hat Y},{\hat H})$ that form a closed 
algebra under
commutation relations. This defines a doubly deformed universal
enveloping algebra ${\cal U}_{q,h}(sl(2))$.

\section {Representations of ${\cal U}_{q,h}(sl(2))$}

The action of the generators of the algebra ${\cal U}_q(sl(2))$ on the
standard basis $\{|jm \rangle|(2j+1) \in  N, -j \leq m \leq j\}$ are
given by
\begin{eqnarray}
&&{\cal J}_{\pm} |jm\rangle= ([j \mp m] [j \pm m+1])^{1/2} |j m \pm 1>,\\
&&q^{\pm {\cal J}_0} |jm\rangle= q^{\pm m} |jm\rangle.
\end{eqnarray}
Repeated actions of ${\cal J}_+$ on the basis states may be expressed as
\begin{equation}
{\cal J}_+^p |jm \rangle = \left(\frac{[j-m] ! [j+m+p] !}{[j+m] ! [j-m-p]
    !}\right)^{1/2} |j\ m+p\rangle
\end{equation}
where $[n] ! = \prod_{k=1}^n [k]$.
Using (41) and the map (10), the action of
the operator ${\hat X}$ on the basis states is obtained as
\begin{equation}
{\hat X} | jm\rangle= \sum_{k \geq 0} \left(\frac{h}{2}\right)^{2k}
\frac{P_k(\xi)}{[2k+1]} \left(\frac{[j-m]!
    [j+m+2k+1]!}{[j+m]![j-m-2k-1]!}\right)^{1/2} |j\ m+2k+1\rangle 
\end{equation}
To obtain the action of ${\hat Y}$ on the basis states, we use the identity
(9) and reexpress ${\hat Y}$ in the normal ordered form :
\begin{eqnarray}
&&{\hat Y} = \sum_{k,l\geq 0} \left(- \frac{h^2}{4}\right)^{k+l}
\biggl({1/2 \atop k}\biggr)  \biggl({1/2 \atop l}\biggr) \\
&&\left({\cal J}_- {\cal J}_+^{2(k+l)} + [2k] {\cal J}_+^{2k+2l-1} 
[2{\cal J}_0 + 2k + 4l -1]\right).
\end{eqnarray}

Using (41), it may be readily obtained
\begin{eqnarray}
&&{\hat Y} |jm\rangle= \sum_{k,l\geq 0} \left(- \frac{h^2}{4}\right)^{k+l}
\biggl({1/2 \atop k}\biggr)  \biggl({1/2 \atop l}\biggr)   \\
&&\left\{\left(\frac{[]j+m+2k+2l][j-m-2k-2l+1][j-m]![j+m+2k+2l]!}{[j+m]!
[j-m-2k-2l]!}\right)^{1/2}\right.\\
&&+ [2k]
[2k+4l+2m-1]\left(\left.\frac{[j-m]![j+m+2k+2l-1]!}{[j+m]![j-m-2k-2l+1]!}
\right)^{1/2}\right\}
 |j\ m+2k+2l-1\rangle \qquad\qquad  
\end{eqnarray}
The Casimir operator for the algebra ${\cal U}_q(sl(2))$
\begin{eqnarray}
&&C_q = {\cal J}_+{\cal J}_- + [{\cal J}_0] [{\cal J}_0-1]\\
&&\phantom{C_q }  = {\cal J}_-{\cal J}_+ + [{\cal J}_0] [{\cal J}_0+1]
\end{eqnarray}
may be expressed in terms of the generators $({\hat X},{\hat Y},
{\hat H})$, leading to a
$q$-deformation of (6). But the eigenvalues will, of
course, be given by
\begin{eqnarray}
C_q |jm\rangle= [j] [j+1] |\ jm\rangle.
\end{eqnarray}
This completes our construction of the irreducible representations of
${\cal U}_{q,h}(sl(2))$ for generic values of $q$.

\section{Conclusion}

We have discussed here nonlinear mappings starting from the algebras
$sl(2)$ and ${\cal U}_q(sl(2))$, leading to the algebras ${\cal U}_h(sl(2))$ 
and
${\cal U}_{q,h}(sl(2))$ respectively. Elsewhere [13] the relation between
$sl(2)$ and ${\cal U}_q(sl(2))$ was also exhibited as the limiting case of a
class of nonlinear maps (where other references can be found). The maps 
leading from $sl(2)$ to ${\cal U}_q(sl(2))$ and ${\cal U}_h(sl(2))$ are 
complementary in the following sense :

\smallskip 

{\bf (i).} In the case of $q$-deformation the nonlinearity enters through the
diagonalizable generator $J_0$.

\smallskip

{\bf (ii).} For the $h$-deformation it enters via the nilpotent
(non-diagonalizable) generator $J_+$.

\smallskip

Finally both cases are combined in the passage to ${\cal U}_{q,h}(sl(2))$. In
terms of the nonlinear maps, here we have followed the path
$$
{\cal U}(sl(2))\longrightarrow {\cal U}_q(sl(2))\longrightarrow 
{\cal U}_{q,h}(sl(2)).
$$
It is also possible to envisage the alternate route
$$
{\cal U}(sl(2)) \longrightarrow {\cal U}_{h}(sl(2)) 
\longrightarrow {\cal U}_{q,h}(sl(2)) .
$$

This is less simple. One has to implement (in the rhs of (10) and (11)) the map
[10, 16] relating (1) and (7) and then the inverse of the map (2). Thus, for 
example,
\begin{eqnarray}
{\hat X}= {2\over h} \sum_{n\geq 0} \alpha_n \biggl(\hbox{tanh}
({h X\over 2})\biggl({[{\hat J}+{1\over 2}]^2-[H+{1\over 2}]^2 \over 
({\hat J}+{1\over 2})^2-(H+{1\over 2})^2 }\biggr)^{1/2}\biggr)^{2n+1}
\end{eqnarray}
where the operator ${\hat J}$ is here given in terms of $C(X,Y,H)$ of (6) as
\begin{eqnarray}
{\hat J}({\hat J}+1)=C.
\end{eqnarray}
One obtains ${\hat Y}$ analogously.  
The following pattern of possibilities emerges:
$$
\diagram{
{\cal U}(sl(2)) & \hfl{\leftarrow}{} & {\cal U}_{h}(sl(2)) \cr
\vfl{\uparrow}{} && \vfl{\uparrow}{} \cr
{\cal U}_{q}(sl(2)) 
& \hfl{\leftarrow}{} & 
{\cal U}_{q,h}(sl(2))  \cr}
$$

For following the small arrows, it is more simple at each stage to take 
the limits
($q \rightarrow 1$ and/or $h \rightarrow 0$), instead of inverting the
map. The problem of explicit construction of irreducible
representations is quite efficiently solved at each higher stage via our maps.

An induced coalgebraic structure of the algebra ${\cal U}_{q,h}(sl(2))$ is
immediately generated by the map (10,11,12). The induced coproduct for the 
generators read (in terms of coproducts (8))
\begin{eqnarray}
\label{formule 4.1}
&& \Delta({\hat X}) = \frac{2}{h} \sum_{n\geq 0} \alpha_n \left(\frac{h}{2}
  \Delta ({\cal J}_+)\right)^{2n+1}\\
&&\Delta({\hat Y}) = \left(1 -
  \left(\frac{h}{2}\Delta({\cal J}_+)\right)^2\right)^{1/2} \Delta({\cal J}_-)\left(1 -
  \left(\frac{h}{2}\Delta({\cal J}_+)\right)^2\right)^{1/2}\\
&&\Delta({\hat H}) = \Delta({\cal J}_0)
\end{eqnarray}
The full induced Hopf structure may be obtained. 

An alternative induced Hopf structure is provided by the map indicated by 
(51) and the Hopf structure of ${\cal U}_h(sl(2))$ [2]. The 
noncocommutative 
coproducts for ${\cal U}_h(sl(2))$ are
\begin{eqnarray}
&& \Delta( X) =  X\otimes 1 +1\otimes X\\
&&\Delta( Y) = Y \otimes e^{h X}+ e^{-h X}\otimes Y\\
&&\Delta( H) = H \otimes e^{h X}+ e^{-h X}\otimes H
\end{eqnarray}
Thus 
\begin{eqnarray}
&& \Delta'( {\hat X}) = {2 \over h} \sum_{n \geq 0} \alpha_{n}
\biggl(\hbox{tanh}({h \over 2}  
\Delta( X))\biggl({[\Delta({\hat J})+{1\over 2}(1\otimes 1)]^2-[\Delta(H)+
{1\over 2}(1\otimes 1)]^2 
\over (\Delta({\hat J})+{1\over 2}(1\otimes 1))^2-(\Delta(H)+{1\over 2}
(1\otimes 1))^2 }
\biggr)^{1/2}\biggr)^{2n+1} 
\end{eqnarray}
The corresponding coproducts $\Delta'({\hat Y})$ and $\Delta'({\hat H})$
(and also the counits and antipodes for all) can thus be written directly. 
Appropriate (explicitly given) inverse maps implemented in the rhs for 
$\Delta$ and $\Delta'$, complete de picture.   
The ${\cal R}$-matrices of ${\cal U}_h(sl(2))$ (see [10] and references 
cited there) can be implemented here through $\Delta'$ (whereas the well-known 
${\cal R}$ matrices for the standard $q$-deformation correspond to $\Delta$). 
The ``non-classical'' automorphism of ${\cal U}_h(sl(2))$ can also be 
implemented via (51).

New features and possibilities arising for $q$ a root of unity remain to 
be studied. Another possibility is worth mentioning. It is known
that at the level of one parameter deformations, the $q$ and $h$
deformations are essentially distinct. The $h$-deformation is a
non-invertible singular limit \cite{AKS95} of the $q$-deformation. Can the
$(q,h)$-deformation may be obtained as a singular limit of a
$(q,q')$-deformation, where $q$ and $q'$ play comparable roles in the
algebra. In this context, we mention that a class of double
deformations was presented in [13] with ${\cal U}_q(sl(2))$ as the starting
point.

Recently we have studied elsewhere [17,18] other two parametric deformations
of $sl(2)$ involving elliptic functions. The discussion in [18] concerning 
possibilities of applications has evident parallel features in the present 
case. Instead of more or less reproducing it we just mention that the {\em 
difference} of the role of the parameters (($h$,$k$) in [18] and ($q$,$h$) 
here) can be of particular interest.


\end{document}